\documentstyle[12pt,amssymb]{article}
\newcommand{\beq}{\begin{equation}}
\def\rank{\hbox{\rm{rank}}}
\newcommand{\eeq}{\end{equation}}

\newtheorem{theorem}{Theorem}[section]

\newtheorem{conjecture}[theorem]{Conjecture}
\newtheorem{proposition}[theorem]{Proposition}
\newtheorem{corollary}[theorem]{Corollary}
\newtheorem{lemma}[theorem]{Lemma}

\begin{document}

\title{A MATRIX AND ITS INVERSE: REVISITING MINIMAL RANK COMPLETIONS}
\author{Hugo J. Woerdeman\footnote{Supported in part by National Science Foundation grant
DMS-0500678}}

\maketitle

\begin{center}
Department of Mathematics\\
Drexel University\\
         3141 Chestnut Street\\
         Philadelphia, PA 19104\\
hugo@math.drexel.edu
\end{center}

\begin{abstract}
We revisit a formula that connects the minimal ranks of triangular
parts of a matrix and its inverse and relate the result to
structured rank matrices.  We also address the generic minimal
rank problem.
\end{abstract}


\bigskip

\noindent{\bf MR Classification:} 15A09, 15A15, 65F05.

\section{Introduction}

In this paper we revisit the following result from \cite{W1}:

\it Let $[(T_{ij})_{i,j=1}^n]^{-1} = (S_{ij})_{i,j=1}^n$ be block
matrices with sizes that are compatible for multiplication. Other
than the full matrix (which is of size $N$, say), none of the
blocks need to be square. Then
\begin{equation}\label{mr}
 \hbox{\rm min \ rank } \pmatrix{ T_{11} &  ? & \cdots & ? \cr T_{21} & T_{22} & \cdots & ? \cr
\vdots &  & \ddots & \vdots \cr T_{n1} & T_{n2} & \cdots & T_{nn}
} + \hbox{\rm min \ rank } \pmatrix{ ? & ? & \cdots & ? \cr S_{21}
& ? & \cdots & ? \cr \vdots & \ddots & \ddots & \vdots \cr S_{n1}
& \cdots & S_{n,n-1} & ?  } =N. \end{equation} \rm With the recent
interest in numerical algorithms that make effective use of
matrices with certain rank structures (see, e.g., \cite{BGP},
\cite{VVM}, \cite{Ty}, \cite{EGO}, \cite{FG}, and references
therein), it seems appropriate to revisit this formula that
captures many of the rank considerations that go into these
algorithms. The nullity theorem due to \cite{Gu} is a particular
case. The papers \cite{SN} and \cite{VV} show the recent interest
in the nullity theorem. It is our hope that this general formula
(\ref{mr}) enhances the insight in rank structured matrices.

In addition, in Section 3 we will address the so-called "generic
minimal rank problem". This problem was introduced by Professors
Gilbert Strang and David Ingerman.

\section{Minimal ranks of matrices and their inverses}

Let us recall the notion of partial matrices and their minimal
rank. Let $\mathbb {F}$ be a field and let $n, m, \nu_{1}, \ldots,
\nu_{n}, \mu_{1}, \ldots, \mu_{m}$ be nonnegative integers. The
{\it pattern} of specified entries in a partial matrix will be
described by a set ${J} \subset$ $ \big\{ {1, \ldots, n}\big\}
\times \big\{{1, \ldots, m} \big\}$. Let now $A_{ij}, (i, j) \in
{J}$, be given matrices with entries in $\mathbb {F}$ of size
$\nu_{i} \times \mu_{j}$. We will allow $\nu_{i}$ and $\mu_{j}$ to
equal $0$.  The collection of matrices $\mathcal {A} = \big\{$
$A_{ij};$ $(i,j) \in {J} \big\}$ is called a {\it partial block
matrix with the pattern} ${J}.$ A block matrix $B =
(B_{ij})^{n}_{i=1{,}}{}^{m}_{j=1}$ with $B_{ij}\in \mathbb {F}$
${\!\!\!\:}^{\nu_{i} \times \mu_{j}}$ is called a {\it completion}
of $\mathcal A$ if $B_{ij} = A_{ij}, (i, j)\in {J}.$ The {\it
minimal rank} of $A$ (notations: min rank$(A)$) is defined by min
rank$(A) = {\rm min} \big\{{\rm rank} \,B : B \, {\rm is \,a
\,completion \,of} \,A \big\}.$

A completion of $\mathcal {A}$ with rank min rank$(A)$ is called a
minimal rank completion of $A$.  When all the blocks are of size
$1 \times 1$ (i.e., $\nu_{i} = \mu_{j}=1$ for all $i$ and $j$), we
will simply talk about a {\it partial matrix}.  Clearly, any block
matrix as above may be viewed as a partial matrix of size $N
\times M$ as well, where $N = \nu_{1} + \ldots + \nu_{n}, M =
\mu_{1} + \ldots + \mu_{m}.$ It will be convenient to represent
partial block matrices in matrix format.  As usual a question mark
will represent an unknown block.  For instance, $\mathcal {A} =$
$\big\{A_{ij}:$ $1 \leq j \leq i \leq n \big\}$ will be
represented as
$$\mathcal {A} = \pmatrix{A_{11} & ? & \ldots & ? \cr \vdots & \ddots & \ddots & \vdots \cr
\vdots & & \ddots & ? \cr A_{n1} & \ldots & \ldots & A_{nn}} $$

The formula that connects the minimal ranks of triangular parts of
a matrix and its inverse is the following. The result appeared
originally in \cite{W1} (see also \cite{W2} and Chapter 5 of
\cite{W}).

\begin{theorem} \cite{W1}\label{inv}
Let $T = (T_{ij})^{n}_{i, j=1}$ be an invertible block matrix with
$T_{ij}$ of size $\nu_{i} \times \mu_{j},$ where $\nu_{i} \geq 0,$
$\mu_{j} \geq 0$ and $N = \nu_{i} + \ldots + \nu_{n} = \mu_{i} +
\ldots + \mu_{n}.$  Put $T^{-1} = (S_{ij})^{n}_{1, j=1}$ where
$S_{ij}$ is of size $\mu_{i} \times \nu_{j}.$ Then $$ \hbox{\rm
min \ rank } \pmatrix{ T_{11} &  ? & \cdots & ? \cr T_{21} &
T_{22} & \cdots & ? \cr \vdots &  & \ddots & \vdots \cr T_{n1} &
T_{n2} & \cdots & T_{nn} } + \hbox{\rm min \ rank } \pmatrix{ ? &
? & \cdots & ? \cr S_{21} & ? & \cdots & ? \cr \vdots & \ddots &
\ddots & \vdots \cr S_{n1} & \cdots & S_{n,n-1} & ?  } =N. $$
\end{theorem}

As we will see, one easily deduces from Theorem \ref{inv} that the
inverse of an upper Hessenberg matrix has the lower triangular
part of a rank 1 matrix. The strength of Theorem \ref{inv} lies in
that one easily deduces a multitude of such results from it.

From the same paper \cite{W1} we would also like to recall the
following result.

\begin{theorem}\cite{W1}\label{2}
The partial matrix ${\mathcal T} = \big\{ T_{ij}:1 \leq j \leq i
\leq n \big\}$ has minimal rank $${\rm min \, rank \,}{\mathcal T}
= \sum_{i=1}^{n} {\rm rank} \pmatrix{T_{i1} & \ldots & T_{ii} \cr
\vdots & & \vdots \cr T_{n1} & \ldots & T_{ni}} - \sum_{i=1}^{n-1}
{\rm rank} \pmatrix{T_{i+1, 1} & \ldots & T_{i+1, i} \cr \vdots &
& \vdots \cr T_{n1} & \ldots & T_{ni}}.$$
\end{theorem}

After the $n=2$ case of Theorem \ref{2} is obtained it is
straightforward to prove the general case by introduction.  For
the $2 \times 2$ case of Theorem \ref{2} one needs to observe that
the minimal rank of
$$\pmatrix{T_{11} & ? \cr T_{21} & T_{22}}
$$
will at least be the rank of $\pmatrix{T_{11} \cr T_{21}}$ plus
the minimal number of columns in $T_{22}$ that together with the
columns of $T_{21}$ span the column space of $\pmatrix{T_{21} &
T_{22}}.$ Once such a minimal set of columns in $T_{22}$ has been
identified, put any numbers on top of these columns.  Now any
other columns in $T_{22}$ can be completed to be a linear
combination of fully completed columns. Doing this leads to a
completion of rank
$${\rm rank} \pmatrix{T_{11} \cr T_{21}} + {\rm rank} \pmatrix{T_{21} & T_{22}} - {\rm
rank} T_{21},$$ yielding the $n=2$ case of Theorem \ref{2}.

The proof of Theorem \ref{inv}, which can be found in \cite{W1} is
easily derived from Theorem \ref{2} and the nullity theorem, which
we recall now.

\begin{theorem}\label{nullity}\cite{Gu}
Consider $$\pmatrix{A & B \cr C & D}^{-1} = \pmatrix{P & Q \cr R &
S}.$$ Then ${\rm dim \, ker \,} C = {\rm dim \, ker\,} R . $
\end{theorem}

{\bf Proof.} Since $CP=-DR$ , $P[\ker R]\subseteq\ker C$.
Likewise, since $RA=-SC$, we get $A[\ker C]\subseteq \ker R$.
Consequently,
$$AP[\ker R]\subseteq A[\ker C]\subseteq \ker R.$$
Since $AP+BR=I$ , $AP[\ker R]=\ker R$, thus
$$A[\ker C]=\ker R.$$
This yields $\dim\ker C\ge \dim\ker R$.  By reversing the roles of
$C$ and $R$ one obtains also that $\dim\ker R\ge \dim\ker C$. This
gives $\dim\ker R=\dim\ker C$, yielding the lemma.\quad $\square$

The nullity theorem is in fact the $n=2$ case of Theorem 1.
Indeed, if $$T^{-1} = \pmatrix{T_{11} & T_{12} \cr T_{21} &
T_{22}}^{-1} = \pmatrix{S_{11} & S_{12} \cr S_{21} & S_{22}},
$$ we get from Theorem 1 that \begin{equation}\label{star}
{\rm rank} \pmatrix{T_{11} \cr T_{21}} + {\rm rank}
\pmatrix{T_{21} & T_{22}} - {\rm rank} T_{21} + {\rm rank} S_{21}
= N.\end{equation} As $T$ is invertible we have that
$\pmatrix{T_{11} \cr T_{21}}$ and $\pmatrix{T_{21} & T_{22}}$ are
full rank, so (\ref{star}) gives
$$\mu_{1} + \nu_{2} - {\rm rank} \,T_{21} + {\rm rank} \,S_{21} = \mu_{1} + \mu_{2} = \nu_{1} + \nu_{2},$$
and thus $$\nu_{2} - {\rm rank} \,T_{21} = \mu_{2} - {\rm rank}
\,S_{21},$$ which is exactly Theorem 3.

To make the connection with some of the results in the literature
we need the following proposition.

\begin{proposition}\label{ext}
Let $\mathcal {T} =$ $\big\{t_{ij}:$ $1\leq j \leq i \leq n\big\}$
be a scalar valued partial matrix.  Then ${\rm min \, rank \,}
(\mathcal{T}) =$ $n$ if and only if $t_{ii} \neq 0, i = 1, \ldots,
n,$ and $t_{ij}=0$ for $i>j.$
\end{proposition}

{\bf Proof.} The "if" part is immediate.  For the only if part
write
\begin{equation}\label{dag} {\rm min\, rank} {\mathcal{T}}={\rm
rank} \pmatrix{t_{11} \cr \vdots \cr t_{n1}} +
\sum_{i=2}^{n}s_{i},\end{equation} where $$s_{i} = {\rm rank}
\pmatrix{t_{i1} & \ldots & t_{ii} \cr \vdots & & \vdots \cr t_{n1}
& \ldots & t_{ni}} - {\rm rank} \pmatrix{t_{i1} & \ldots &
t_{i,i-1} \cr \vdots \cr t_{n1} & \ldots & t_{n, i-1}}. $$ All the
terms in (\ref{dag}) are at the most 1, and as there are exactly
$n$ terms they need to all be equal to 1 for min
rank$(\mathcal{T})=$ $n$ to be satisfied.  But then $s_{n}=1$
implies $t_{n1}=\ldots=t_{n,n-1}=0$ and $t_{nn} \neq 0.$
Inductively, one can then show that $s_{k}=1$ implies
$t_{k1}=\ldots=t_{k,k-1}=0$ and $t_{kk} \neq 0,k={n-1}, \ldots,
2.$ Finally the first column of $\mathcal{T}$ needs to have rank
1.  As $t_{ij} = 0, j=2, \ldots, n,$ was already established we
get that $t_{11}\neq 0.$ This proves the result. $\square$

We now easily obtain the following corollary, due to Asplund
\cite{As}.

\begin{corollary} \cite{As}
Let $p \ge 0$ and $A=(a_{ij})^{N}_{i,j=1}$ be an $N \times N$
scalar matrix with inverse $B=(b_{ij})^{N}_{i,j=1}.$ Then $a_{ij}
= 0$ for all $i$ and $j$ with $j > i+p$, and $a_{ij} \neq 0$,
$j=i+p$ if and only if there exist a $N\times p$ matrix $F$ and a
$p\times N$ matrix G so that $b_{ij} = (FG)_{ij}$, $ i < j+p$.
\end{corollary}

{\bf Proof}. Let $(S_{ij})_{i,j=N-p+1} = A$, where $S_{i1}$ is of
size $1 \times p$, $i=1,\ldots , n-p$, $S_{N-p+1, 1}$ has size $p
\times p$, $S_{N-p+1, j}$ has size $p \times 1$, $j=2, \ldots ,
N-p+1$, and all the other $S_{ij}$ are $1\times 1$. Let
$B=(T_{ij})_{i,j=1}^{N-p+1}$ be partitioned accordingly. Then, it
follows from (\ref{mr}) that
$$ \hbox{\rm min \ rank } \pmatrix{ ? & ? & \cdots & ? \cr S_{21} &
? & \cdots & ? \cr \vdots & \ddots & \ddots & \vdots \cr S_{n1} &
\cdots & S_{n,n-1} & ?  } = N-p $$ if and only if
$$ \hbox{\rm min \ rank } \pmatrix{ T_{11} &  ? & \cdots & ? \cr
T_{21} & T_{22} & \cdots & ? \cr \vdots &  & \ddots & \vdots \cr
T_{n1} & T_{n2} & \cdots & T_{nn} } = p . $$ Using Proposition
\ref{ext} the result now follows. $\square$

In a similar way it is easy to deduce results by \cite{BF},
\cite{R}, \cite{RBRF}, \cite{RBFF1}, \cite{RRB} and \cite{El} from
Theorem \ref{inv}. For instance, if $T_{ij}$ and $S_{ij}$ are
scalars, and $T_{21}, \ldots, T_{n,n-1} \neq 0$ and $T_{ij} =0$
for $i>j+1$, then the left hand term in (\ref{mr}) is $n-1$. Since
$N=n$, we get that the lower triangular partial matrix
$(S_{ij})_{i\ge j}$ has minimal rank 1. Thus one easily obtains
that $S_{ij} = u_i v_j, i\ge j$, where $u_1, \ldots , u_n, v_1,
\ldots , v_n$ are scalars. Examples like this show that Theorem
\ref{inv} is useful in the contexts of semi-separability and
quasi-separability (see, e.g., \cite{VVGM} and \cite{EG} for an
overview of these notions). We hope that the simplicity of formula
(\ref{mr}) will help in the further development of these notions.

\section{The generic minimal rank completion problem}
Recently D. Ingerman and G. Strang posed the following problem.
Suppose that a partial matrix (over some field $\mathbb F$) has
the property that all of its fully specified submatrices are of
full rank and so that every $k \times k$ partial submatrix has at
most $(2k-r)r$ entries specified. Is it true that one can always
complete to a matrix of rank $\le r$? The count of $(2k-r)r$
specified entries comes from the consideration that if $r$ columns
and $r$ rows in a $k\times k$ submatrix are specified, one can
complete this submatrix to a rank $r$ one (due to the fact that
the submatrix in the overlap of the $r$ columns and the $r$ rows
has full rank). However, as soon as one adds one specified entry
to these $r$ columns and $r$ rows, immediately a $(r+1)\times
(r+1)$ submatrix is specified, and the minimal rank will be at
least $r+1$.

Ingerman and Strang showed that the above statement is correct for
$r=1$. However, the following example shows that in general it is
not correct for $r \ge 2$.

{\bf Example.} Consider the matrix
$$ A:= \pmatrix{ 6 & 3 & x & 1 \cr 3 & 1 & 1 & y \cr z & 1 & 2 & 3 \cr
1 & w & 1 & 1 } , $$ where $x,y,z$ and $w$ are the unknowns. Note
that this partial matrix satisfies the requirements stated in the
first paragraph. Furthermore, suppose that rank$A= 2$. Then we
have that
$$ \pmatrix{ 6 & 3 \cr 3 & 1 }  - \pmatrix{ x & 1 \cr
1 & y } \pmatrix{ 2 & 3 \cr 1 & 1 }^{-1} \pmatrix{ z & 1 \cr 1 & w
} = 0 , $$ and since the rank of the first term is 2, the second
term must also have rank 2. Thus, we have that $xy\neq1$ and $zw
\neq 1$. Next, we also have that
$$ \pmatrix{ z & 1 \cr 1 & w } - \pmatrix{ 2 & 3 \cr 1 & 1 } \pmatrix{ x & 1 \cr
1 & y }^{-1} \pmatrix{ 6 & 3 \cr 3 & 1 } = 0 . $$ Multiplying on
both sides with $xy-1$, the off-diagonal entries yield the
following equations
$$ xy - 6y -3x +10 =0, \ \ \ xy -6y -3x +8 = 0 . $$
These are not simultaneously solvable (as long as we are in a
field where $8 \neq 10$). It should be noted that this is a
counterexample for any field in which $6 \neq 9$, $6\neq 1$, $3
\neq 1$, $9 \neq 1$ (so that we have full rank specified
submatrices) and $8 \neq 10$. As an aside, we note that for some
of the small fields it may impossible to fulfill the nondegeneracy
requirement on the data. E.g., when ${\mathbb F} = \{ 0,1 \}$ a
$2\times2$ matrix can only be nonsingular if zeroes are allowed in
the matrix.

It should be noted that if one associates the bipartite graph with
the partial matrix (see, e.g., \cite{CJRW}) one obtains a minimal
eight cycle. Consequently, the bipartite graph is not bipartite
chordal as bipartite chordality requires by definition the absence
of minimal cycles of length 6 or greater. Notice that in the $r=1$
case the condition on the density of the specified entries
prevents the existence of minimal cycles of length 6 or more. We
now arrive at the following conjecture.

\begin{conjecture}\label{conj} Consider a partial matrix for which the bipartite graph is
bipartite chordal. Suppose furthermore that any fully specified
submatrix has full rank and that any $k \times k$ submatrix has at
most $(2k-r)r$ entries specified. Then there exists a completion
of rank $r$.
\end{conjecture}

We can prove the conjecture for the subclass of banded patterns
(cf. \cite{W3}).

\begin{theorem} \label{th:2.1}
Consider a partial matrix with a banded pattern (as defined in
\cite{W3}). Suppose furthermore that any fully specified submatrix
has full rank and that any $k \times k$ submatrix has at most
$(2k-r)r$ entries specified. Then there exists a completion of
rank $r$.
\end{theorem}

{\em Proof.} By Theorem 1.1 in \cite{W3} it suffices to show that
for every triangular subpattern (for the definition, see
\cite{W3}) we have that the minimal rank is $\le r$. But a
triangular subpattern can always embedded in a pattern that
corresponds to $r$ rows and columns specified (due to the
condition that in any $k \times k$ submatrix has at most $(2k-r)r$
entries are specified). But then the result follows. $\Box$

\bigskip

Observe that the proof shows that if the bipartite chordal minimal
rank conjecture in \cite{CJRW} (see also Chapter 5 in \cite{W}) is
true, then the above conjecture is true as well. The techniques
developed in \cite{BB} and/or \cite{JM} may be helpful in proving
the conjecture above.

\end{document}